\documentclass[12pt,a4paper]{article}
\usepackage{amsmath}
\usepackage{amssymb}
\usepackage{t1enc}
\usepackage[latin1]{inputenc}
\usepackage[english]{babel}
\usepackage{amsfonts}
\usepackage{latexsym}
\usepackage{bm}
\usepackage{setspace}
\usepackage[toc,page]{appendix}
\usepackage{graphicx}
\usepackage{enumerate}
\usepackage[numbers,sort&compress]{natbib}
\usepackage{tikz}

\usepackage{theoremref}
\usepackage{lineno}

\usepackage{mathtools}

\newtheorem{theorem}{Theorem}%[section]

\newtheorem{lemma}{Lemma}

\newtheorem{construction}{Construction}

\allowdisplaybreaks

\begin{document}
\title{Proof of a conjecture on isolation of graphs dominated by a vertex 
\medskip\medskip
}

\author{Peter Borg \\[2mm]
\normalsize Department of Mathematics \\
\normalsize Faculty of Science \\
\normalsize University of Malta\\
\normalsize Malta\\
\normalsize \texttt{peter.borg@um.edu.mt}
}

\date{}%{\today} 
\maketitle

\begin{abstract}
A copy of a graph $F$ is called an $F$-copy. For any graph $G$, the $F$-isolation number of $G$, denoted by $\iota(G,F)$, is the size of a smallest subset $D$ of the vertex set of $G$ such that the closed neighbourhood $N[D]$ of $D$ in $G$ intersects the vertex sets of the $F$-copies contained by $G$ (equivalently, $G-N[D]$ contains no $F$-copy). Thus, $\iota(G,K_1)$ is the domination number $\gamma(G)$ of $G$, and $\iota(G,K_2)$ is the vertex-edge domination number of $G$. We prove that if $F$ is a $k$-edge graph, $\gamma(F) = 1$ (that is, $F$ has a vertex that is adjacent to all the other vertices of $F$), and $G$ is a connected $m$-edge graph, then $\iota(G,F) \leq \big\lfloor \frac{m+1}{k+2} \big\rfloor$ unless $G$ is an $F$-copy or $F$ is a $3$-path and $G$ is a $6$-cycle. This was recently posed as a conjecture by Zhang and Wu, who settled the extreme case where $F$ is a star. The result for the other extreme case where $F$ is a clique had been obtained by Fenech, Kaemawichanurat and the present author. The bound is attainable for any $m \geq 0$ unless $1 \leq m = k \leq 2$. New ideas, including deletion methods and divisibility considerations, are introduced in the proof of the conjecture. 
\end{abstract}

\section{Introduction} \label{Introsection}
Unless stated otherwise, we use capital letters such as $X$ to denote sets or graphs, and small letters such as $x$ to denote non-negative integers or elements of a set. The set of positive integers is denoted by $\mathbb{N}$. For $n \geq 1$, $[n]$ denotes the set $\{1, \dots, n\}$ (that is, $\{i \in \mathbb{N} \colon i \leq n\}$). We take $[0]$ to be the empty set $\emptyset$. Arbitrary sets are taken to be finite. For a set $X$, the set of $2$-element subsets of $X$ is denoted by ${X \choose 2}$ (that is, ${X \choose 2} = \{ \{x,y \} \colon x,y \in X, x \neq y \}$). For standard terminology in graph theory, we refer the reader to \cite{West}. Most of the terminology used here is defined in \cite{Borg}. 

Every graph $G$ is taken to be \emph{simple}, that is, its vertex set $V(G)$ and edge set $E(G)$ satisfy $E(G) \subseteq {V(G) \choose 2}$. We may represent an edge $\{v,w\}$ by $vw$. We call $G$ an \emph{$n$-vertex graph} if $|V(G)| = n$. We call $G$ an \emph{$m$-edge graph} if $|E(G)| = m$. For $v \in V(G)$, $N_{G}(v)$ denotes the set of neighbours of $v$ in $G$, $N_{G}[v]$ denotes the closed neighbourhood $N_{G}(v) \cup \{ v \}$ of $v$, and $d_{G}(v)$ denotes the degree $|N_{G}(v)|$ of $v$. For $X \subseteq V(G)$, $N_G[X]$ denotes the closed neighbourhood $\bigcup_{v \in X} N_G[v]$ of $X$, $G[X]$ denotes the subgraph of $G$ induced by $X$, and $G - X$ denotes the graph obtained by deleting the vertices in $X$ from $G$. Thus, $G[X] = (X,E(G) \cap {X \choose 2})$ and $G - X = G[V(G) \backslash X]$. Where no confusion arises, the subscript $G$ may be omitted; for example, $N_G(v)$ may be abbreviated to $N(v)$. With a slight abuse of notation, for $Y \subseteq E(G)$, $G - Y$ denotes the graph obtained by removing the edges in $Y$ from $G$, that is, $G - Y = (V(G), E(G) \backslash Y)$. For $x \in V(G) \cup E(G)$, we may abbreviate $G - \{x\}$ to $G - x$. 

A \emph{component of $G$} is a maximal connected subgraph of $G$. Clearly, the components of $G$ are pairwise vertex-disjoint, and their union is $G$. 

Consider two graphs $G$ and $H$. If $G$ is a copy of $H$, then we write $G \simeq H$ and we say that $G$ is an \emph{$H$-copy}. If $H$ is a subgraph of $G$, then we say that \emph{$G$ contains $H$}. %We say that $G$ is \emph{$H$-free} if $G$ contains no copy of $H$.  
  
For $n \geq 1$, the graphs $([n], {[n] \choose 2})$, $([n], \{\{1, i\} \colon i \in [n] \backslash \{1\}\})$ and $([n], \{\{i,i+1\} \colon i \in [n-1]\})$ are denoted by $K_n$, $K_{1,n-1}$ and $P_n$, respectively. For $n \geq 3$, $C_n$ denotes the graph $([n], \{\{1,2\}, \{2,3\}, \dots, \{n-1,n\}, \{n,1\}\})$. A $K_n$-copy is called an \emph{$n$-clique} or a \emph{complete graph}, a $K_{1,n}$-copy is called a \emph{star}, a $P_n$-copy is called an \emph{$n$-path} or simply a \emph{path}, and a $C_n$-copy is called an \emph{$n$-cycle} or simply a \emph{cycle}. %A $3$-cycle is a $3$-clique and is also called a \emph{triangle}. 

If $D \subseteq V(G) = N[D]$, then $D$ is called a \emph{dominating set of $G$}. The size of a smallest dominating set of $G$ is called the \emph{domination number of $G$} and is denoted by $\gamma(G)$. If $\mathcal{F}$ is a set of graphs and $F$ is a copy of a graph in $\mathcal{F}$, then we call $F$ an \emph{$\mathcal{F}$-graph}. If $D \subseteq V(G)$ such that $G-N[D]$ contains no $\mathcal{F}$-graph, then $D$ is called an \emph{$\mathcal{F}$-isolating set of $G$}. Note that $D$ is an $\mathcal{F}$-isolating set of $G$ if and only if $N[D]$ intersects the vertex sets of the $\mathcal{F}$-graphs contained by $G$. Let $\iota(G, \mathcal{F})$ denote the size of a smallest $\mathcal{F}$-isolating set of $G$. If $\mathcal{F} = \{F\}$, then we may replace $\mathcal{F}$ in these defined terms and notation by $F$. Clearly, $D$ is a $K_1$-isolating set of $G$ if and only if $D$ is a dominating set of $G$. Thus, $\gamma(G) = \iota(G, K_1)$.

The study of isolating sets was introduced by Caro and Hansberg~\cite{CaHa17}. It is a natural generalization of the study of dominating sets \cite{C, CH, HHS, HHS2, HL, HL2}. One of the earliest results in this field is the upper bound $n/2$ of Ore \cite{Ore} on the domination number of any connected $n$-vertex graph $G \not\simeq K_1$ (see \cite{HHS}). While deleting the closed neighbourhood of a dominating set yields the graph with no vertices, deleting the closed neighbourhood of a $K_2$-isolating set yields a graph with no edges. In the literature, a $K_2$-isolating set is also called a \emph{vertex-edge dominating set}. Consider any connected $n$-vertex graph $G$. Caro and Hansberg~\cite{CaHa17} proved that $\iota(G, K_2) \leq n/3$ unless $G \simeq K_2$ or $G \simeq C_5$. This was independently proved by \.{Z}yli\'{n}ski \cite{Z} and solved a problem in \cite{BCHH} (see also \cite{LMS, BG}). Let $\mathcal{C}$ be the set of cycles. Solving one of the problems posed by Caro and Hansberg \cite{CaHa17}, the present author~\cite{Borg} proved that 
\begin{equation} \iota(G,\mathcal{C}) \leq \frac{n}{4} \label{Borgcyclebound}
\end{equation} 
unless $G \simeq K_3$, and that the bound is sharp. A special case of this result is that $\iota(G,K_3) \leq n/4$ unless $G \simeq K_3$. Solving another problem posed in \cite{CaHa17}, Fenech, Kaemawichanurat and the present author~\cite{BFK} proved that 
\begin{equation} \iota(G, K_k) \leq \frac{n}{k+1} \label{BFKbound}
\end{equation} 
unless $G \simeq K_k$ or $k = 2$ and $G \simeq C_5$, and that the bound is sharp. The result of Ore and the result of Caro and Hansberg and of \.{Z}yli\'{n}ski are the cases $k = 1$ and $k = 2$, respectively. In \cite{CaHa17}, it was also shown that $\iota(G,K_{1,k}) \leq \frac{n}{k+1}$. For $k \geq 1$, let $\mathcal{F}_{0,k} = \{K_{1,k}\}$, let $\mathcal{F}_{1,k}$ be the set of regular graphs of degree at least $k-1$, let $\mathcal{F}_{2,k}$ be the set of graphs whose chromatic number is at least $k$, and let $\mathcal{F}_{3,k} = \mathcal{F}_{0,k} \cup \mathcal{F}_{1,k} \cup \mathcal{F}_{2,k}$ (see \cite{Borgrsc, Borgrc2}). In \cite{Borgrsc}, the present author proved that for each $i \in \{0, 1, 2, 3\}$, 
\begin{equation} \iota(G,\mathcal{F}_{i,k}) \leq \frac{n}{k+1}  \label{Borggenbound}
\end{equation}
unless $i \neq 0$ and either $G \simeq K_k$ or $k = 2$ and $G \simeq C_5$, and that if $i \neq 0$, then the bound is sharp. This generalizes all the results above as $\mathcal{C} \subseteq \mathcal{F}_{1,3}$ and $K_k \in \mathcal{F}_{1,k} \cap \mathcal{F}_{2,k}$. Other isolation bounds of this kind in terms of $n$ are given in \cite{BBS, Borgcon, CX, Yan, ZW2, ZW3}. It is worth mentioning that domination and isolation have been particularly investigated for maximal outerplanar graphs \cite{BK, BK2, CaWa13, CaHa17, Ch75, DoHaJo16, DoHaJo17, HeKa18, LeZuZy17, Li16, MaTa96, KaJi, To13}, mostly due to connections with Chv\'{a}tal's Art Gallery Theorem \cite{Ch75}. As in the development of domination, isolation is expanding in various directions; for example, \emph{total} isolation and an isolation game have been treated in \cite{BGH} and \cite{BDJKR}, respectively. 

Consider any connected $m$-edge graph $G$. Fenech, Kaemawichanurat and the present author~\cite{BFK2} also proved that, analogously to (\ref{BFKbound}), 
\begin{equation} \iota(G, K_k) \leq \frac{m+1}{{k \choose 2} + 2} \label{BFKbound2}
\end{equation} 
unless $G \simeq K_k$. Let $\mathcal{F}_{1,k}$ and $\mathcal{F}_{2,k}$ be as above, and let $\mathcal{F}_{3,k}$ now be $\mathcal{F}_{1,k} \cup \mathcal{F}_{2,k}$. Recently, the present author \cite{Borgrc2} proved that, analogously to (\ref{Borggenbound}), for each $i \in \{1, 2, 3\}$, 
\begin{equation} \iota(G, \mathcal{F}_{i,k}) \leq \frac{m+1}{{k \choose 2} + 2} \label{Borggenbound2}
\end{equation} 
unless $G \simeq K_k$, and that the bound is sharp. This generalizes (\ref{BFKbound2}) and immediately yields $\iota(G,\mathcal{C}) \leq \frac{m+1}{5}$ if $G \not\simeq K_3$. In this paper, we prove a conjecture of Zhang and Wu \cite{ZW} that generalizes (\ref{BFKbound2}) in another direction. Before presenting the result,  we construct a graph that attains the conjectured bound. The following is a generalization of \cite[Construction 1.2]{BFK2} and a slight variation of the construction of $B_{n,F}$ in \cite{Borg}. 

\begin{construction} \label{const2} \emph{Consider any $m, k \in \{0\} \cup \mathbb{N}$ and any connected $k$-edge graph $F$, where $F \simeq K_1$ if $k = 0$ (that is, $V(F) \neq \emptyset$). By the division algorithm, there exist $q, r \in \{0\} \cup \mathbb{N}$ such that $m+1 = q(k+2) + r$ and $0 \leq r \leq k+1$. Let $Q_{m,k}$ be a set of size $q$. If $q \geq 1$, then let $v_1, \dots, v_q$ be the elements of $Q_{m,k}$, let $F_1, \dots, F_q$ be copies of $F$ such that the $q+1$ sets $V(F_1), \dots, V(F_q)$ and $Q_{m,k}$ are pairwise disjoint, and for each $i \in [q]$, let $w_i \in V(F_i)$, and let $G_i$ be the graph with $V(G_i) = \{v_i\} \cup V(F_i)$ and $E(G_i) = \{v_iw_i\} \cup E(F_i)$. If either $q = 0$, $T$ is the null graph $(\emptyset, \emptyset)$, and $G$ is a connected $m$-edge graph $T'$, or $q \geq 1$, $T$ is a tree with vertex set $Q_{m,k}$ (so $|E(T)| = q-1$), $T'$ is a connected $r$-edge graph with $V(T') \cap \bigcup_{i=1}^q V(G_i) = \{v_q\}$, and $G$ is a graph with $V(G) = V(T') \cup \bigcup_{i=1}^q V(G_i)$ and $E(G) = E(T) \cup E(T') \cup \bigcup_{i=1}^q E(G_i)$, then we say that $G$ is an \emph{$(m,F)$-special graph} with \emph{quotient graph $T$} and \emph{remainder graph $T'$}, and for each $i \in [q]$, we call $G_i$ an \emph{$F$-constituent of $G$}, and we call $v_i$ the \emph{$F$-connection of $G_i$ in $G$}. We say that an $(m,F)$-special graph is \emph{pure} if its remainder graph has no edges (\cite[Figure~1]{BFK2} is an illustration of a pure $(71, K_5)$-special graph). Clearly, an $(m,F)$-special graph is a connected $m$-edge graph.}
\end{construction}

If $F$ and $G$ are graphs such that either $G \simeq F$ or $F \simeq K_{1,2}$ ($\simeq P_3$) and $G \simeq C_6$, then we say that $(G,F)$ is \emph{special}. In the next section, we prove the following result.

\begin{theorem} \label{mainresult} If $F$ is a $k$-edge graph with $\gamma(F) = 1$,  $G$ is a connected $m$-edge graph, and $(G,F)$ is not special, then 
\[\iota(G, F) \leq \left \lfloor \frac{m+1}{k+2} \right \rfloor.\]
Moreover, equality holds if $G$ is an $(m,F)$-special graph.
\end{theorem}
This proves Conjecture~4.4 of the recent paper \cite{ZW}, in which Zhang and Wu treated the case where $F$ is a star. Note that $\gamma(F) = 1$ means that $F$ has a vertex $v$ that is adjacent to all the other vertices of $F$, so $N_F[v] = V(F)$ and $F$ is connected. Also note that the bound in Theorem~\ref{mainresult} is attained for any $m \geq 0$ unless $1 \leq m = k \leq 2$. This follows from the second part of the theorem because clearly it is only when $1 \leq m = k \leq 2$ that every $(m,F)$-special graph $G$ is such that $(G,F)$ is special. Note that the bound is attained if $m = 0$ and $G = (\emptyset, \emptyset)$.

New ideas, including deletion methods and divisibility considerations, are introduced in the proof of Theorem~\ref{mainresult} out of necessity. They are mostly concentrated in the argument for Case~1.2 of the proof. These ideas promise to be of further use in the rapidly expanding study of graph isolation.

\section{Proof of Theorem~\ref{mainresult}} \label{Proofsection}

We start the proof of Theorem~\ref{mainresult} with two lemmas from \cite{Borg}.

\begin{lemma} \label{lemma}
If $G$ is a graph, $\mathcal{F}$ is a set of graphs, $X \subseteq V(G)$, and $Y \subseteq N[X]$, then $\iota(G, \mathcal{F}) \leq |X| + \iota(G-Y, \mathcal{F})$. 
\end{lemma}
\textbf{Proof.} Let $D$ be an $\mathcal{F}$-isolating set of $G-Y$ of size $\iota(G-Y, \mathcal{F})$. Clearly, $V(F) \cap Y \neq \emptyset$ for each $\mathcal{F}$-graph $F$ that is a subgraph of $G$ and not a subgraph of $G-Y$. Since $Y \subseteq N[X]$, $X \cup D$ is an $\mathcal{F}$-isolating set of $G$. The result follows.~\hfill{$\Box$}

\begin{lemma} \label{lemmacomp}
If $G_1, \dots, G_r$ are the distinct components of a graph $G$, and $\mathcal{F}$ is a set of connected graphs, then $\iota(G,\mathcal{F}) = \sum_{i=1}^r \iota(G_i,\mathcal{F})$.
\end{lemma}
\textbf{Proof.} For each $i \in [r]$, let $D_i$ be a smallest $\mathcal{F}$-isolating set of $G_i$. Consider any $\mathcal{F}$-graph $F$ contained by $G$. Since $F$ is connected, there exists some $j \in [r]$ such that $G_j$ contains $F$, so $N[D_j] \cap V(F) \neq \emptyset$. Thus, $\bigcup_{i = 1}^r D_i$ is an $\mathcal{F}$-isolating set of $G$, and hence $\iota(G, \mathcal{F}) \leq \sum_{i = 1}^r |D_i| = \sum_{i = 1}^r \iota(G_i, \mathcal{F})$. Let $D$ be a smallest $\mathcal{F}$-isolating set of $G$. For each $i \in [r]$, $D \cap V(G_i)$ is an $\mathcal{F}$-isolating set of $G_i$. We have $\sum_{i = 1}^r \iota(G_i, \mathcal{F}) \leq \sum_{i = 1}^r |D \cap V(G_i)| = |D| = \iota(G, \mathcal{F})$. The result follows.~\hfill{$\Box$}
\\

For a vertex $v$ of a graph $G$, let $E_G(v)$ denote the set $\{vw \colon w \in N_G(v)\}$. For $X, Y \subseteq V(G)$, let $E_G(X,Y)$ denote the set $\{xy \in E(G) \colon x \in X, \, y \in Y\}$.
\\

\noindent
\textbf{Proof of Theorem~\ref{mainresult}.} We first prove the second part of the theorem. Thus, suppose that $G$ is an $(m,F)$-special graph with exactly $q$ $F$-constituents as in Construction~\ref{const2}, and that $(G, F)$ is not special. Then, $|E(G)| = m$ and $q = \left \lfloor \frac{m+1}{k+2} \right \rfloor$. If $q = 0$, then $m \leq k$, and hence, since $G \not\simeq F$, $\iota(G,F) = 0 = q$. Suppose $q \geq 1$. Then, $\{v_1, \dots, v_q\}$ is an $F$-isolating set of $G$, so $\iota(G, F) \leq q$. If $D$ is an $F$-isolating set of $G$, then, since $G_1 - v_1, \dots, G_q - v_q$ are copies of $F$, $D \cap V(G_i) \neq \emptyset$ for each $i \in [q]$. Therefore, $\iota(G, F) = q$.  

We now prove that the bound in the theorem holds. Thus, let $G$ be a connected $m$-edge graph such that $(G, F)$ is not special, and let $n = |V(G)|$ and $\ell = |V(F)|$. Since $\iota(G, F)$ is an integer, it suffices to prove that $\iota(G, F) \leq \frac{m+1}{k+2}$. We use induction on $n$. If $k$ is $0$ or $1$, then $F \simeq K_1$ or $F \simeq K_2$, and hence the result is given by (\ref{BFKbound2}). Suppose $k \geq 2$. Then, $\ell \geq 3$.  Let $\mathcal{S}$ be the set of $F$-copies contained by $G$. If $\mathcal{S} = \emptyset$, then $\iota(G,F) = 0 \leq \frac{m+1}{k+2}$. Suppose $\mathcal{S} \neq \emptyset$. Then, $\iota(G, F) \geq 1$ and $n \geq \ell$. Since $\gamma(F) = 1$, for each $S \in \mathcal{S}$, $V(S) = N_S[v_S]$ for some $v_S \in V(S)$. Let $U = \{u \in V(G) \colon V(S) = N_S[u] \mbox{ for some } S \in \mathcal{S}\}$. Let $v \in U$ such that $d_G(u) \leq d_G(v)$ for each $u \in U$. For some $F_1 \in \mathcal{S}$, $V(F_1) = N_{F_1}[v] \subseteq N_G[v]$. Thus, $d(v) \geq \ell - 1 \geq 2$. Since $G$ is connected and contains $F_1$, and $(G, F)$ is not special (so $G \neq F_1$), $m \geq |E(F_1)| + 1 = k+1$. If $V(G) = N[v]$, then $\iota(G, F) = 1 \leq \frac{m+1}{k+2}$. Suppose $V(G) \neq N[v]$. Let $G' = G-N[v]$ and $n' = |V(G')|$. Then, $V(G') \neq \emptyset$.

Let $\mathcal{H}$ be the set of components of $G'$. For any $H \in \mathcal{H}$ and any $x \in N(v)$ such that $xy_{x,H} \in E(G)$ for some $y_{x,H} \in V(H)$, we say that $H$ is \emph{linked to $x$} and that $x$ is \emph{linked to $H$}. Since $G$ is connected, for each $H \in \mathcal{H}$, $x_Hy_H \in E(G)$ for some $x_H \in N(v)$ and some $y_H \in V(H)$. We have 
\begin{equation} E(F_1) \subseteq E(G[N[v]]), \quad \{x_Hy_H \colon H \in \mathcal{H}\} \subseteq E_G(N(v), V(G')), \label{edgesubsets1}
\end{equation}
\begin{equation} m = |E(G[N[v]])| + |E_G(N(v), V(G'))| + \sum_{H \in \mathcal{H}} |E(H)| \geq k + \sum_{H \in \mathcal{H}} |E(H) \cup \{x_Hy_H\}|. \label{sizeineq_1}
\end{equation} 
Let $\mathcal{H}' = \{H \in \mathcal{H} \colon (H,F) \mbox{ is special}\}$. By the induction hypothesis, $\iota(H, F) \leq \frac{|E(H)|+1}{k+2}$ for each $H \in \mathcal{H} \backslash \mathcal{H}'$. 
\\

\noindent
\textbf{Case 1:} \emph{$\mathcal{H}' = \emptyset$.} By Lemma~\ref{lemma} (with $X = \{v\}$ and $Y = N[v]$) and Lemma~\ref{lemmacomp},
\begin{align}
\iota(G, F) &\leq 1 + \iota(G', F) = 1 + \sum_{H \in \mathcal{H}} \iota(H, F) \label{iotaineq_1} \\
&\leq \frac{k+2}{k+2} + \sum_{H \in \mathcal{H}} \frac{|E(H)|+1}{k+2}. \nonumber
\end{align}
Thus, if $m \geq k + 1 + \sum_{H \in \mathcal{H}} (|E(H)|+1)$, then $\iota(G, F) \leq \frac{m+1}{k+2}$. Suppose $m < k + 1 + \sum_{H \in \mathcal{H}} (|E(H)|+1)$. Then, by (\ref{sizeineq_1}), $m = k + \sum_{H \in \mathcal{H}} (|E(H)|+1)$ and 
\begin{equation} E(G) = E(F_1) \cup \bigcup_{H \in \mathcal{H}} (E(H) \cup \{x_Hy_H\}). \label{extremeH}
\end{equation} 
We have $(k+2) \iota(H, F) \leq |E(H)| + 1$ for each $H \in \mathcal{H}$.\medskip

\noindent
\textbf{Case 1.1:} \emph{$(k+2) \iota(I, F) \leq |E(I)|$ for some $I \in \mathcal{H}$.} Then, 
\begin{align} m + 1 &\geq k + 2 + (k+2) \iota(I, F) + \sum_{H \in \mathcal{H} \backslash \{I\}} (k+2) \iota(H,F) \nonumber \\
&= (k+2) \left( 1 + \sum_{H \in \mathcal{H}} \iota(H,F) \right) \geq (k+2) \iota(G, F) \quad \mbox{(by (\ref{iotaineq_1}))}, \nonumber
\end{align}
so $\iota(G, F) \leq \frac{m+1}{k+2}$.\medskip

\noindent
\textbf{Case 1.2:} \emph{$(k+2) \iota(H, F) = |E(H)| + 1$ for each $H \in \mathcal{H}$.} Here we introduce the idea of considering $I - y_I$ for some member $I$ of $\mathcal{H} \backslash \mathcal{H}'$ (which, in this case, is $\mathcal{H}$ as we are in Case 1), and the idea of considering $G - (\{v\} \cup V(I))$ and $G - V(I)$. Let $I' = I - y_I$. Let $\mathcal{J}$ be the set of components of $I'$. Since $I$ is connected, for each $J \in \mathcal{J}$, $y_Iz_J \in E(G)$ for some $z_J \in V(J)$. Let $\mathcal{J}' = \{J \in \mathcal{J} \colon (J,F) \mbox{ is special}\}$. For each $J \in \mathcal{J} \backslash \mathcal{J}'$, let $D_J$ be an $F$-isolating set of $J$ of size $\iota(J,F)$. By the induction hypothesis, $|D_J| \leq \frac{|E(J) \cup \{y_Iz_J\}|}{k+2}$ for each $J \in \mathcal{J} \backslash \mathcal{J}'$.

Suppose $\mathcal{J'} = \emptyset$. Let $G'' = G - (\{v\} \cup V(I))$. Let $L_1, \dots, L_s$ be the distinct components of $G''$. Consider any $i \in [s]$. By (\ref{extremeH}), there exist some $\mathcal{H}_i \subseteq \mathcal{H}$ and $A_i \subseteq E(F_1) \backslash E_G(v)$ (each of $\mathcal{H}_i$ and $A_i$ is possibly empty) such that 
\[E(L_i) = A_i \cup \bigcup_{H \in \mathcal{H}_i} (E(H) \cup \{x_Hy_H\}).\] 
Thus, 
\begin{equation} |E(L_i)| = |A_i| + \sum_{H \in \mathcal{H}_i} (|E(H)| + 1) = |A_i| + \sum_{H \in \mathcal{H}_i} (k+2) \iota(H,F). \label{divisibility}
\end{equation}
We now introduce divisibility considerations. Since $A_i \subseteq E(F_1) \backslash E_G(v)$ and $d(v) \geq 2$, 
$|A_i| \leq k - d(v) \leq k - 2$. Thus, by (\ref{divisibility}), $(L_i,F)$ is not special (otherwise, $|E(L_i)| = k$ or $(k, |E(L_i)|) = (2, 6)$). By the induction hypothesis, 
\[\iota(L_i,F) \leq \frac{|E(L_i)| + 1}{k+2} = \frac{|A_i| + 1}{k+2} + \sum_{H \in \mathcal{H}_i} \iota(H,F) \leq \frac{k-1}{k+2} + \sum_{H \in \mathcal{H}_i} \iota(H,F).\]
Since $\iota(L_i,F)$ is an integer, we obtain $\iota(L_i,F) \leq \sum_{H \in \mathcal{H}_i} \iota(H,F)$. Now, clearly, $\bigcup_{i = 1}^s \mathcal{H}_i = \mathcal{H} \backslash \{I\}$ (and $\mathcal{H}_1, \dots, \mathcal{H}_s$ are pairwise disjoint), so we have
\[\sum_{i=1}^s \iota(L_i,F) \leq \sum_{i=1}^s \sum_{H \in \mathcal{H}_i} \iota(H,F) = \sum_{H \in \mathcal{H} \backslash \{I\}} \iota(H,F) = \sum_{H \in \mathcal{H} \backslash \{I\}} \frac{|E(H) \cup \{x_Hy_H\}|}{k+2}.\]
By (\ref{extremeH}), the components of $G - \{v, y_I\}$ are $L_1, \dots, L_s$ and the members of $\mathcal{J}$. Since $v, y_I \in N[x_I]$, Lemmas~\ref{lemma} and \ref{lemmacomp} give us
\begin{align} \iota(G,F) &\leq 1 + \iota(G - \{v, y_I\}) = 1 + \sum_{i=1}^s \iota(L_i,F) + \sum_{J \in \mathcal{J}} \iota(J,F) \nonumber \\
&\leq \frac{|E(F_1) \cup \{x_Iy_I\}| + 1}{k+2} + \sum_{H \in \mathcal{H} \backslash \{I\}} \frac{|E(H) \cup \{x_Hy_H\}|}{k+2} + \sum_{J \in \mathcal{J}} \frac{|E(J) \cup \{y_Iz_J\}|}{k+2} \nonumber \\
&\leq \frac{m+1}{k+2}. \nonumber
\end{align}

Now suppose $\mathcal{J'} \neq \emptyset$. If $J \in \mathcal{J}'$ and $J \simeq F$, then let $D_J = \emptyset$. If $J \in \mathcal{J}'$ and $J \simeq C_6$ (so $F \simeq K_{1,2}$), then let $z_J'$ be the vertex of $J$ such that $V(J) = N_J[\{z_J, z_J'\}]$, and let $D_J = \{z_J'\}$. Let $G^* = G - V(I)$. Then, $G^*$ is connected and contains $F_1$. Suppose $G^* \simeq C_6$. Then, $E(G^*) = \{vw_1, w_1w_2, w_2w_3, w_3w_4, w_4w_5, w_5v\}$ for some distinct $w_1, \dots, w_5 \in V(G) \backslash \{v\}$. Let $H^* = (\{w_2, w_3, w_4\}, \{w_2w_3, w_3w_4\})$. Since $V(F_1) \subseteq N[v]$, (\ref{extremeH}) gives us that $V(F_1) = N[v] = \{v, w_1, w_5\}$, $H^* \in \mathcal{H}$, and $H^*$ is linked to $x_{H^*}$ only, which contradicts $w_1w_2, w_4w_5 \in E(G^*)$. Thus, $G^* \not\simeq C_6$. If $\mathcal{H} = \{I\}$, then $G^* = F_1$ and we let $D = \emptyset$. If $\mathcal{H} \neq \{I\}$, then $G^* \neq F_1$ and we let $D$ be an $F$-isolating set of $G^*$ of size $\iota(G^*, F)$. By the induction hypothesis, $|D| \leq \frac{|E(G^*)|+1}{k+2}$. Clearly, $D \cup \{y_I\} \cup \bigcup_{J \in \mathcal{J}} D_J$ is an $F$-isolating set of $G$, so
\begin{align} \iota(G, F) &\leq |D| + 1 + \sum_{J \in \mathcal{J}'} |D_J| + \sum_{J \in \mathcal{J} \backslash \mathcal{J}'} |D_J| \nonumber \\
&\leq \frac{|E(G^*)|+1}{k+2} + \frac{|\{x_Iy_I\} \cup \bigcup_{J \in \mathcal{J'}} (E(J) \cup \{y_Iz_J\})|}{k+2} + \sum_{J \in \mathcal{J} \backslash \mathcal{J}'} \frac{|E(J) \cup \{y_Iz_J\}|}{k+2} \nonumber \\
&\leq \frac{m+1}{k+2}. \nonumber
\end{align}
\vskip 2mm

\noindent
\textbf{Case 2: $\mathcal{H}' \neq \emptyset$.} For each $x \in N(v)$, let $\mathcal{H}'_x = \{H \in \mathcal{H}' \colon H \mbox{ is linked to } x\}$ and $\mathcal{H}_x^* = \{H \in \mathcal{H} \backslash \mathcal{H}' \colon H \mbox{ is linked to $x$ only}\}$. For each $H \in \mathcal{H} \backslash \mathcal{H}'$, let $D_H$ be an $F$-isolating set of $H$ of size $\iota(H,F)$.

Suppose $H \simeq C_6$ for some $H \in \mathcal{H}'$. Then, $F \simeq K_{1,2}$ and $H = (\{y_1, \dots, y_6\}$, $\{y_1y_2, \dots, y_5y_6, y_6y_1\})$ for some distinct $y_1, \dots, y_6 \in V(G)$ with $y_1 = y_H$. Let $G^* = G - N_H[y_4]$ and $A^* = \{y_2y_3, y_3y_4,y_4y_5, y_5y_6\}$. Then, $G^*$ is connected and $A^* \subseteq E(G) \backslash E(G^*)$. We have $x_H, y_2, y_6 \in N_{G^*}(y_1)$, so $y_1 \in U$ (as $F \simeq K_{1,2}$), and hence $d_G(v) \geq d_G(y_1) \geq 3$. Thus, since $d_{G^*}(v) = d_G(v)$, $G^*$ is neither a copy of $F$ nor a $6$-cycle. Let $D^*$ be a smallest $F$-isolating set of $G^*$. By the induction hypothesis, $|D^*| \leq \frac{|E(G^*)|+1}{k+2} = \frac{|E(G^*)|+1}{4}$. By Lemma~\ref{lemma} (with $X = \{y_4\}$ and $Y = N_H[y_4]$),
\begin{align} \iota(G,F) &\leq 1 + |D^*| \leq \frac{|E(G^*)| + 5}{4} = \frac{|E(G^*)| + |A^*| + 1}{k+2} \leq \frac{m+1}{k+2}. \nonumber
\end{align} 

Now suppose $H \not\simeq C_6$ for each $H \in \mathcal{H}'$. Then, each member of $\mathcal{H}'$ is a copy of $F$.\medskip

\noindent
\textbf{Case 2.1:} \emph{$|\mathcal{H}'_x| \geq 2$ for some $x \in N(v)$.} Let $X = \{x_H \colon H \in \mathcal{H}' \backslash \mathcal{H}'_x\}$. We have $x \notin X \subset N(v)$, so $d(v) \geq 1 + |X|$. Let $D = \{v, x\} \cup X \cup \left( \bigcup_{H \in \mathcal{H} \backslash \mathcal{H}'} D_{H} \right)$. We have $V(G) = N[v] \cup \bigcup_{H \in \mathcal{H}} V(H)$, $y_{x,H} \in N[x]$ for each $H \in \mathcal{H}'_x$, and $y_{H} \in N[x_H]$ for each $H \in \mathcal{H}' \backslash \mathcal{H}'_x$, so $D$ is an $F$-isolating set of $G$. Since $\iota(G,F) \leq |D|$ and
\begin{align} m + 1 &\geq 1 + d(v) + \sum_{H \in \mathcal{H}'_x} |E(H) \cup \{x y_{x,H}\}| + \sum_{H \in \mathcal{H} \backslash \mathcal{H}_x'} |E(H) \cup \{x_H y_{H}\}| \nonumber \\
&\geq 2 + |X| + (k+1)|\mathcal{H}_x'| + (k+1)|\mathcal{H}' \backslash \mathcal{H}_x'| + \sum_{H \in \mathcal{H} \backslash \mathcal{H}'} (k+2)|D_H| \nonumber \\
&\geq 2 + |X| + 2(k+1) + (k+1)|X| + \sum_{H \in \mathcal{H} \backslash \mathcal{H}'} (k+2)|D_H| = (k+2)|D|, \nonumber
\end{align}
$\iota(G,F) \leq \frac{m+1}{k+2}$.\medskip

\noindent
\textbf{Case 2.2:}
\begin{equation} |\mathcal{H}'_x| \leq 1 \mbox{ \emph{for each} } x \in N(v). \label{k=3Hx<2} 
\end{equation} 
Let $H \in \mathcal{H}'$. Let $x = x_H$ and $y = y_{H}$.\medskip

\noindent
\textbf{Case 2.2.1:} \emph{$H$ is linked to $x$ only.} Let $X = \{x\} \cup V(H)$. Then, $G - X$ has a component $G_v^*$ such that $N[v] \backslash \{x\} \subseteq V(G_v^*)$, and the other components of $G - X$ are the members of $\mathcal{H}_{x}^*$. Let $D^*$ be an $F$-isolating set of $G_v^*$ of size $\iota(G_v^*,F)$, and let $D = \{x\} \cup D^* \cup \bigcup_{I \in \mathcal{H}_x^*} D_I$. Then, $D$ is an $F$-isolating set of $G$. Since 
\begin{equation} E(G) \supseteq \{vx, xy\} \cup E(H) \cup E(G_v^*) \cup \bigcup_{I \in \mathcal{H}_x^*} (E(I) \cup \{xy_{x,I}\}), \nonumber  %\label{starpoint1}
\end{equation}
\begin{equation} m+1 \geq 3 + k + |E(G_v^*)| + \sum_{I \in \mathcal{H}_x^*} (|E(I)| + 1) \geq 3 + k + |E(G_v^*)| + \sum_{I \in \mathcal{H}_x^*} (k+2)|D_I|. \label{starpoint2}
\end{equation}

Suppose that $(G_v^*, F)$ is not special. By the induction hypothesis, $|D^*| \leq \frac{|E(G_v^*)|+1}{k+2}$. By (\ref{starpoint2}), 
\[m+1 \geq 3 + k + (k+2)|D^*| - 1 + \sum_{I \in \mathcal{H}_x^*} (k+2)|D_I| = (k+2)|D| \geq (k+2)\iota(G,F),\]
so $\iota(G,F) \leq \frac{m+1}{k+2}$. 

Now suppose that $(G_v^*,F)$ is special. Suppose first that $G_v^* \simeq F$. Then, $\{x\} \cup \bigcup_{I \in \mathcal{H}_x^*} D_I$ is an $F$-isolating set of $G$, and hence $\iota(G, F) \leq 1 + \sum_{I \in \mathcal{H}_x^*} |D_I|$. By (\ref{starpoint2}),
\[m + 1 \geq 3 + 2k + \sum_{I \in \mathcal{H}_x^*} (k+2)|D_I| > (k+2)\iota(G, F),\] 
so $\iota(G, F) < \frac{m+1}{k+2}$. Now suppose $G_v^* \not\simeq F$. Then, $F \simeq K_{1,2}$ and $G_v^* \simeq C_6$. Thus, $G_v^* = (\{v, x_1, x_2, w_1, w_2, w_3\}$, $\{vx_1, x_1w_1, w_1w_2, w_2w_3, w_3x_2, x_2v\})$, where $\{x_1, x_2\} = N(v) \backslash \{x\}$ and $w_1, w_2, w_3$ are the distinct elements of $V(G) \backslash (N[v] \cup X \cup \bigcup_{I \in \mathcal{H}_x^*} V(I))$. Since $\{x, w_2\} \cup \bigcup_{I \in \mathcal{H}_x^*} D_I$ is an $F$-isolating set of $G$, $\iota(G, F) \leq 2 + \sum_{I \in \mathcal{H}_x^*} |D_I|$. By (\ref{starpoint2}), $m + 1 \geq 9 + k + \sum_{I \in \mathcal{H}_x^*} (k+2)|D_I|$. Since $k = 2$, $\iota(G, F) < \frac{m+1}{k+2}$.\medskip

\noindent
\textbf{Case 2.2.2:} \emph{$H$ is linked to some $x' \in N(v) \backslash \{x\}$.} Then, $x'y' \in E(G)$ for some $y' \in V(H)$. Let $I = G - V(H)$. Then, $I$ is connected. Since $H \simeq F$, $V(H) \subseteq N[w]$ for some $w \in V(H)$. Let $A = E_G(N(v), V(H))$. Then, $xy, x'y' \in A$.

Suppose that $(I, F)$ is not special. By Lemma~\ref{lemma} (with $X = \{w\}$ and $Y = V(H)$) and the induction hypothesis, 
\[\iota(G,F) \leq 1 + \iota(I,F) \leq \frac{|E(H) \cup A|}{k+2} + \frac{|E(I)|+1}{k+2} \leq \frac{m+1}{k+2}.\]

Now suppose that $(I, F)$ is special. Suppose $I \simeq C_6$. Then, $F \simeq K_{1,2}$. Thus, $k = 2$, $F_1 \simeq K_{1,2}$ and $H \simeq K_{1,2}$. Since $V(F_1) \subseteq N[v]$, it follows that $N(v) = \{x, x'\}$ and $E(I) = \{vx, xz_1, z_1z_2, z_2z_3, z_3x', x'v\}$ with $\{z_1, z_2, z_3\} = V(I) \backslash N[v]$. Since $F \simeq K_{1,2}$ and $v, z_1, y \in N(x)$, we have $x \in U$ and $d(x) \geq 3 > d(v)$, contradicting the choice of $v$. Therefore, $I \not\simeq C_6$. Thus, $I \simeq F$, and hence $I = F_1$. We have $E(G) = E(H) \cup E(I) \cup A$ and $V(H) \cap V(I) = \emptyset$, so
\[m = |E(H)| + |E(I)| + |A| = 2k + |A|.\] 
Since $V(F_1) \subseteq N[v]$ and $V(H) \subseteq N[w]$, we have $V(G) = N[\{v, w\}]$, so $\iota(G, F) \leq 2 = \frac{m - |A| + 4}{k+2}$. Thus, if $|A| \geq 3$, then $\iota(G, F) \leq \frac{m+1}{k+2}$. Suppose $|A| \leq 2$. Since $xy, x'y' \in A$, $A = \{xy, x'y'\}$. 

Suppose $w \in \{y, y'\}$. Then, $|N(w) \cap \{x, x'\}| \geq 1$. Since $V(H) \subseteq N[w]$, we have $|V(G - N[w])| \leq |V(I)| - 1 = \ell - 1$, so $\iota(G,F) = 1 < \frac{m+1}{k+2}$. 

Now suppose $w \notin \{y, y'\}$. We may assume that 
\begin{equation} d_H(y) \leq d_H(y'). \label{yy'}
\end{equation}
Let $J = G - N[x']$. Let $J_H = J[V(J) \cap V(H)]$ and $J_I = J[V(J) \cap V(I)]$. We have $y' \in V(H) \backslash V(J_H)$, $v, x' \in V(I) \backslash V(J_I)$ and 
\begin{equation} E(J) \subseteq E(J_H) \cup E(J_I) \cup \{xy\}. \label{IH}
\end{equation}
Suppose that $J$ contains an $F$-copy $F_2$. Then, $V(F_2) \subseteq N_J[u]$ for some $u \in V(F_2)$. Since $|V(J_H)| \leq \ell - 1$ and $|V(J_I)| \leq \ell - 2$, it follows by (\ref{IH}) that $u = y$ and $V(F_2) = \{x\} \cup V(H - y')$. By (\ref{IH}), $N_{F_2}(x) = \{u\}$. Thus, since $H \simeq F \simeq F_2$, $d_H(z) = 1$ for some $z \in V(H)$. Since $d_H(w) = \ell - 1 \geq 2$, $z \neq w$. Since $V(H) \subseteq N[w]$, $N_H(z) = \{w\}$. Since $\{x\} \cup V(H - y') = V(F_2) \subseteq N_J[u] = N_J[y] \subseteq \{x\} \cup N_H[y]$, we have $V(H-y') \subseteq N_H[y]$, so $z \in \{y, y'\}$ (otherwise, we obtain $w, y \in N_H(z)$, a contradiction). If $z = y$, then $N_H(y) = \{w\}$. If $z = y'$, then $d_H(y) \leq 1$ by (\ref{yy'}), so again $N_H(y) = \{w\}$. Thus, $V(H) = \{w, y, y'\}$. Since $I \simeq F \simeq H$, $V(I) = \{v, x, x'\}$. Since $x \in V(J)$, we have $xx' \notin E(G)$, so $E(I) = \{vx, vx'\}$. Thus, $F \simeq K_{1,2}$. Since $H \simeq F$, $E(H) = \{wy, wy'\}$. Thus, we have $E(G) = \{vx, xy, yw, wy', y'x', x'v\}$, meaning that $G \simeq C_6$, which is a contradiction as $(G,F)$ is not special.~\hfill{$\Box$}
\\
\\
\noindent
\textbf{Acknowledgements.} The author is grateful to the two anonymous referees for checking the paper carefully and providing constructive remarks.

\end{document}